\documentclass[a4paper,12pt]{article}

\usepackage{makeidx}

\usepackage{latexsym}
\usepackage{amscd}
\usepackage{amsmath}
\usepackage{amssymb}

\usepackage{amsthm}
\usepackage{float}

\usepackage{graphicx}
\usepackage{psfrag}

\usepackage{mymacros}

\newcommand{\wa}{c}

\newcommand{\toricu}{\ensuremath{w}}

\newcommand{\edgea}{\ensuremath{k}}

\newcommand{\Irindex}{\ensuremath{l}}

\newcommand{\Izero}{\ensuremath{I^{(0)}}}
\newcommand{\Ionezero}{\ensuremath{I^{(10)}}}
\newcommand{\Ioneone}{\ensuremath{I^{(11)}}}
\newcommand{\unknownY}{Y}
\newcommand{\llangle}{\langle \langle}
\newcommand{\rrangle}{\rangle \rangle}

\newcommand{\pt}{pt}
\newcommand{\Etilde}{\ensuremath{\tilde{E}}}
\newcommand{\identitymatrix}{\ensuremath{1}}

\newcommand{\spint}{\ensuremath{\scI}}
\newcommand{\docJ}{\ensuremath{D}}
\newcommand{\pair}[2]{\langle #1, #2 \rangle}

\newcommand{\motonot}{t}
\newcommand{\motonoq}{q}
\newcommand{\atonot}{\widetilde{t}}

\newcommand{\blsix}{\ensuremath{Y}}

\newcommand{\rt}{\mathrm{right}}
\newcommand{\lt}{\mathrm{left}}
\newcommand{\rtt}{r}
\newcommand{\ltt}{l}
\newcommand{\formal}{\mathrm{formal}}

\newcommand{\XSigma}{{X_\Sigma}}
\title{Stokes Matrix for the Quantum Cohomology
of Cubic Surfaces}
\author{Kazushi Ueda}
\date{}
\begin{document}

\maketitle

\begin{abstract}
We prove the conjectural relation
between the Stokes matrix
for the quantum cohomology
of $X$
and an exceptional collection
generating $\DbX$
when $X$ is a smooth cubic surface.
The proof is based on a toric degeneration
of a cubic surface,
the Givental's mirror theorem for toric manifolds,
and the Picard-Lefschetz theory.
\end{abstract}

\section{Introduction}
Let $X$ be a smooth cubic surface in $\bP^3$.
It is a rational surface
obtained by blowing up $\bP^2$
at six points.
Let
$$
 \phi : X \rightarrow \bP^2
$$
denote this blow-up.
For a given
$\degreed \in H_2(X;\bZ)$,
let
$$
 n(\degreed) = \int_\degreed c_1 (TX) - 1
$$
be the expected dimension
of the moduli space of
stable maps
of genus zero and degree $\degreed$
and
$$
 N_\degreed = \langle 
        \underbrace{[\pt], \ldots, [\pt]}_{n(\degreed) \ \text{times}})
       \rangle_{0,\degreed}
$$
be the Gromov-Witten invariant,
where $[\pt]$ denotes
the Poincar\'{e} dual
of the homology class
of a point.
When $n(\degreed) > 0$,
$N_\degreed$ is
the number of
nodal rational curves in $X$
of degree $\degreed$
passing through $n(\degreed)$ points
in general position.
Although it is possible to count
$N_{\degreed}$
by elementary method
when the degree $\degreed$ is small,
the difficulty grows rapidly
as $\degreed$ increases.

The theory of Gromov-Witten invariants
gives
a new perspective
to this classical problem
of counting $N_{\degreed}$
in enumerative geometry.
Instead of considering $N_\degreed$'s individually,
one can treat them as a whole by
the following generating function
on the cohomology ring $H^*(X;\bC)$:
\begin{equation} \label{eq:dppotential}
 \Fpotential(t_0, \ldots, t_8)
  = \frac{1}{2} t_0^2 t_8 + t_0 (t_1^2 - t_2^2 - \cdots - t_7^2)
     + \sum_{\degreed \in H_2(X;\bZ)}
         N_{\degreed} q^\degreed
          \frac{t_8^{n(\degreed)}}{n(\degreed)!}. 
\end{equation}
Here,
$t_0, t_1, t_2, \ldots, t_7, t_8$
are coordinates
on the total cohomology ring $H^*(X;\bC)$
of $X$
corresponding to the basis
$[X], [H], [E_1], \ldots, [E_6], [\pt]$
of $H^*(X;\bC)$,
where
$[X]$ denotes
the dual
of the fundamental class of $X$,
$[H]$
that of the hyperplane class, and
$[E_i]$
that of 
the exceptional divisor
for $i=1, \ldots, 6$,
and
$
 q^\degreed = \exp[\langle \degreed, [H] \rangle t_1]
              \exp[\langle \degreed, [E_1] \rangle t_2]
              \cdots
              \exp[\langle \degreed, [E_6] \rangle t_7].
$
Then,
this $\Fpotential$
endows $H^*(X;\bC)$
with the structure of a
{\em Frobenius manifold}.

A Frobenius manifold is a complex manifold
whose holomorphic tangent bundle
has a bilinear form and
an associative commutative product structure
satisfying a series of axioms.
%
In the case of a Frobenius manifold
coming from the Gromov-Witten invariants
of a smooth projective variety,
the bilinear form
is given by
the Poincar\'{e} pairing,
and the structure constant
of the product structure
is given by the third derivative
of the generating function
$\Fpotential$.
This product structure
is called the quantum cohomology ring,
whose associativity gives
a highly non-trivial non-linear
differential equation
for $\Fpotential$,
called the {\em WDVV
(Witten-Dijkgraaf-Verlinde-Verlinde)}
equation,
which sometimes determines
all the genus-zero
Gromov-Witten invariants
from a finite numbers of them
\cite{Kontsevich-Manin_reconstruction}.

An important point
in the theory of Frobenius manifolds
is its connection
%
%
%
with isomonodromic deformations.
Given a Frobenius manifold,
one can construct
an isomonodromic family
of ordinary differential equations
on $\bP^1$
of the following form:
\begin{equation}
 \hbar \frac{\partial \unknownY}{\partial \hbar}
   + \frac{1}{\hbar} \scU(\unknownY)
   - \scV(\unknownY)
  = 0,
   \label{eq:intro_hbar-direction}
\end{equation}
\begin{equation}
  \hbar \frac{\partial \unknownY}{\partial t_\flata}
  = \frac{\partial}{\partial t_\flata} \circ \unknownY, \quad
  \flata=0,\ldots,N-1.
   \label{eq:intro_t-direction}
\end{equation}

Here, $\unknownY$ is an unknown function on $\bP^1$
times the Frobenius manifold
taking value in the tangent bundle of
the Frobenius manifold,
$\hbar$ is the coordinate on $\bP^1$,
$N$ is the dimension of the Frobenius manifold,
$\{t_\flata\}_{\flata=0}^{N-1}$ is
the {\em flat coordinate}
of the Frobenius manifold,
$\circ$ denotes the product on the tangent bundle,
and $\scU$, $\scV$ are certain operators
acting on sections of the tangent bundle.
(\ref{eq:intro_hbar-direction}) is an
ordinary differential equation on $\bP^1$
with a regular singularity at infinity
and an irregular singularity at the origin,
and (\ref{eq:intro_t-direction}) gives
its isomonodromic deformation.
If a point
on the Frobenius manifold
is semisimple,
i.e., if there are no nilpotent elements
in the product structure
on the tangent space
at this point,
one can define the {\em monodromy data}
of (\ref{eq:intro_hbar-direction})
at this point,
consisting of the monodromy matrix at infinity,
the Stokes matrix at the origin,
and the connection matrix between
infinity and the origin.
These data do not depend
on the choice of a semisimple point
because of the isomonodromicity.

The following conjecture,
originally due to Kontsevich,
developed by Zaslow \cite{Zaslow}, and
formulated into the following form
by Dubrovin \cite{Dubrovin_GATFM},
reveals a striking connection
between the Gromov-Witten invariants
and the derived category of coherent sheaves:

\begin{conjecture} \label{conj:stokes}

The quantum cohomology of a smooth projective variety $X$
is semisimple
if and only if
the bounded derived category $\DbX$
of coherent sheaves on $X$
is generated as a triangulated category
by an exceptional collection
$( \scE_i )_{i=1}^N$.
In such a case, the Stokes matrix $S$
for the quantum cohomology of $X$
is given by
\begin{equation} \label{eq:conjecture}
S_{ij} = \sum_k (-1)^k \dim \Ext^k(\scE_i, \scE_j).
\end{equation}
\end{conjecture}

An exceptional collection appearing above
is the following:

\begin{definition} \label{def:exceptional_collection}
\begin{enumerate}
\item 
 An object $\scE$ in a triangulated category
 is exceptional if
  $$
   \Ext^i(\scE,\scE)=\left\{\begin{array}{cl}
      \bC & \mbox{if $i=0$,}\\
       0  & \mbox{otherwise.}\\
   \end{array}\right.
  $$
\item
 An ordered set of objects $(\scE_i)_{i=1}^N$
 in a triangulated category
 is an exceptional collection
 if each $\scE_i$ is exceptional and
 $\Ext^k(\scE_i,\scE_j)=0$ for any $i>j$ and any $k$.

\end{enumerate}

\end{definition}


Conjecture \ref{conj:stokes} was previously known to hold
for projective spaces
\cite{Dubrovin_PT2DTFT},
\cite{Guzzetti}
and Grassmannians
\cite{Ueda_SMQCG}.
The main result in this paper is:

\begin{theorem} \label{th:main_theorem}
 Conjecture \ref{conj:stokes} holds
 for smooth cubic surfaces in $\bP^3$.
\end{theorem}


To prove
Theorem \ref{th:main_theorem},
we first use the deformation invariance
of the Gromov-Witten invariants
to reduce the problem
to the case of
a toric surface
$\blsix$
defined in Section \ref{sc:vanishing_cycles}
obtained
by bringing
some of the blowing-up centers
to infinitely-near points.
Then
we can use
the Givental's mirror theorem
\cite{Givental_MTTCI},
which gives an integral representation
of the fundamental solution
to (\ref{eq:intro_hbar-direction}),
(\ref{eq:intro_t-direction})
for $\blsix$
of the following form:

\begin{equation} \label{eq:int_rep}
 I_\Gamma = \int_{\Gamma}
  \phi(x,y;\hbar) \exp[W(x,y)/\hbar] \frac{dx dy}{x y}
\end{equation}

Here,
$\phi(x,y;\hbar)$
is a certain $H^*(\blsix;\bC)$-valued function
of $x$, $y$, $\hbar$ and
the flat coordinate
$\{ t_\flata \}_{\flata=0}^8$,
and $W(x,y)$ is a Laurent polynomial
of $x$ and $y$
whose coefficients depend on the flat coordinate.
$\Gamma$ runs over
certain cycles in $(\bCx)^2$
depending on $\hbar$
and flat coordinate,
which are descending Morse cycles
for $\Re[W(x,y)/\hbar]$
for a suitable choice of a metric on $(\bCx)^2$.
%

Since the integrand
in (\ref{eq:int_rep})
is single-valued,
the monodromy of the integral
comes solely from the monodromy
of the integration cycle
$\Gamma$,
which is determined by
intersection numbers of vanishing cycles
by the Picard-Lefschetz theory.
This gives the
left-hand side of (\ref{eq:conjecture}).
On the right-hand side,
one can find an exceptional collection
generating the derived category of coherent sheaves
on $\blsix$
by combining the theorems
of Beilinson \cite{Beilinson}
and Kapranov-Vasserot \cite{Kapranov-Vasserot}.
The $\Ext$-groups between them
can be computed explicitly,
and exhibit a complete agreement
with the result
on the left-hand side.


{\bf Acknowledgements}:
We thank
H. Iritani
for patiently explaining
mirror symmetry
for toric manifolds.
We also thank
A. Ishii,
T. Kawai,
H. Kawanoue,
H. Nasu,
K. Saito,
A. Takahashi,
and H. Uehara
for valuable discussions and comments.
The author is supported by
JSPS Fellowships for Young Scientists
No.15-5561.

\section{Frobenius manifold and Stokes matrix}
 \label{sc:Frobenius_manifolds}

We review the relation between Frobenius manifolds
and isomonodromic deformations
in this section.
See, e.g.,
\cite{Dubrovin_G2DTFT},
\cite{Dubrovin_PT2DTFT},
\cite{Hertling_FM},
\cite{Hitchin_FM},
\cite{Manin_FM}
for more details.

\begin{definition}
A Frobenius manifold is a quintuple
$(M,\circ,e,E,g)$
satisfying the following axioms:
\begin{enumerate}
 \item $M$ is a complex manifold.
 \item $g:\scT_M \otimes \scT_M \rightarrow \scO_M$ is an
       $\scO_M$-bilinear form.
 \item $\circ:\scT_M \otimes \scT_M \rightarrow \scT_M$ defines
       an associative commutative
       $\scO_M$-algebra structure on $\scT_M$.
 \item $g(X \circ Y, Z) = g(X,Y\circ Z)$
       as an $\scO_M$-linear map
       $
        \scT_M \otimes \scT_M \otimes \scT_M
         \rightarrow \scO_M
       $
       (the Frobenius property).
 \item $g$ is flat, i.e., the Levi-Citiva connection
       $\nabla$ associated to $g$ is flat. \label{cond:flatness}
 \item The $(3,1)$-tensor $\nabla \circ$ is totally symmetric
       (the potentiality condition).
 \item $e \in \Gamma(M, \scT_M)$ is an identity element for $\circ$.
 \item $e$ is flat, i.e., $\nabla e=0$.
 \item $E \in \Gamma(M, \scT_M)$ satisfies
       \begin{equation} \label{eq:linearity}
	\nabla \nabla E=0, 
       \end{equation}
       \begin{equation} \label{eq:homogeneity_1}
	\Lie_E(\circ)=\circ,
       \end{equation}
       and there exists a complex number $D$ such that
       \begin{equation} \label{eq:homogeneity_2}
	\Lie_E(g)=D g.
       \end{equation}
        \label{cond:euler}
\end{enumerate}
Here, $\scO_M$ is the structure sheaf of $M$,
$\scT_M$ is the tangent sheaf of $M$,
and $\Lie_E$ is the Lie derivative
with respect to the vector field $E$.
We call $g$ the metric, $E$ the Euler vector field,
and $D$ the charge of the Frobenius manifold $M$.
\end{definition}
Let $N$ denote
the dimension of $M$.
Since $g$ is flat,
there exists a local coordinate
$\{t_\flata \}_{\flata=0}^{N-1}$
on $M$,
called the {\em flat coordinate},
satisfying
$$
g(\partial_\flata, \partial_\flatb) = \eta_{\flata \flatb}
$$
which is unique up to affine transformations.
Here, $\eta_{\flata \flatb}$ is a constant matrix
and 
$$
\partial_\flata = \frac{\partial}{\partial t_\flata}.
$$
Since the identity vector field
$e$ is flat,
we can choose
the flat coordinate
$\{ t_\flata \}_{\flata=0}^{N-1}$
so that
$$
 e = \partial_0.
$$
Let $\Fstrconst_{\flata \flatb}^\flatc$ be the structure constant
of the product $\circ$
in the flat coordinate:
$$
 \partial_\flata \circ \partial_\flatb
  = \sum_{\flatc = 0}^{N-1}
     \Fpotential_{\flata \flatb}^\flatc \partial_\flatc.
$$
Then the potentiality condition
is that
$$
\sum_{\flate = 0}^{N-1} \eta_{\flata \flate}
 \partial_\flatb \Fstrconst_{\flatc \flatd}^\flate
$$
is totally symmetric
with respect to $\flata, \flatb, \flatc, \flatd$.
From this condition,
it follows that
locally on $M$,
there exists a holomorphic function
$\Fpotential$
such that
\begin{equation}
\sum_{\flatd=0}^{N-1} \eta_{\flata \flatd} \Fstrconst_{\flatb \flatc}^\flatd
 = \partial_\flata \partial_\flatb \partial_\flatc \Fpotential.
     \label{eq:str_const} 
\end{equation}
This function $\Fpotential$
is called the potential of the Frobenius manifold $M$.
Note that
$$
 \eta_{\flata \flatb}
  = \partial_{0} \partial_{\flata} \partial_{\flatb}
     \Fpotential
$$
follows from the Frobenius property.

The associativity
of $\circ$
gives the following
{\em WDVV (Witten-Dijkgraaf-Verlinde-Verlinde) equation}
\begin{equation} \label{eq:WDVV}
  \sum_{e,f = 0}^{N-1} \frac{\partial^3 \Fpotential}
  {\partial t_\flata  \partial t_\flatb \partial t_e}
  \eta^{e f} \frac{\partial^3 \Fpotential}
  {\partial t_f \partial t_\flatc \partial t_\flatd}
 =  \sum_{e,f = 0}^{N-1} \frac{\partial^3 \Fpotential}
  {\partial t_\flata \partial t_\flatd \partial t_e}
  \eta^{e f} \frac{\partial^3 \Fpotential}
  {\partial t_f \partial t_\flatb \partial t_\flatc},
\end{equation}
which is a non-linear
partial differential equation
for the potential $\Fpotential$.
Here, $\eta^{\flata \flatb}$
is the inverse matrix
of $\eta_{\flata \flatb}$ :
$
\sum_{\flatc=0}^{N-1} \eta^{\flata \flatc} \eta_{\flatc \flatb}
 = \delta^\flata_\flatb.
$


Now let us define
$\scU, \scV \in \Gamma(\scEnd(\scT_M))$
by
\begin{eqnarray}
 \scU(X) &=& E \circ X, \label{eq:U} \\
 \scV(X) &=& \nabla_X E-\frac{D}{2} X \label{eq:V}, 
\end{eqnarray}
and let
$\pi : \bCx \times M \rightarrow M$
be the second projection.
Define a connection
$$
\fscon : \pi^* \scT_M \rightarrow \pi^* \scT_M \otimes
 \Omega_{\bCx \times M}^1
$$
on $\pi^* \scT_M$,
where $\Omega_{\bC^\times \times M}^1$
denotes the cotangent sheaf of
$\bCx \times M$,
by
\begin{eqnarray}
 \fscon_X Y &=& \nabla_X Y
  - \frac{1}{\hbar} X \circ Y,
  \label{eq:X-direction} \\
 \fscon_{\partial_\hbar} Y &=& \partial_\hbar Y
   + \frac{1}{\hbar^2} \scU(Y)
   - \frac{1}{\hbar} \scV(Y),
  \label{eq:hbar-direction}
\end{eqnarray}
for $X$ and $Y$
which are pull-backs
of local sections of $\scT_M$
to $\bCx \times M$
and $\hbar$ is the inhomogeneous coordinate
of $\bP^1$,
and extend it by the Leibniz rule.
Note the $\hbar$ here is the notation by Givental
and corresponds to
$1 / z$ in Dubrovin's papers \cite{Dubrovin_G2DTFT},
\cite{Dubrovin_PT2DTFT}.
$\fscon$ is called
the {\em first structure connection}
of the Frobenius manifold $M$.
The following Theorem \ref{th:fsflat}
is a keystone
in the theory of Frobenius manifolds:
\begin{theorem} \label{th:fsflat}
 $\fscon$ is flat.
\end{theorem}
For the proof,
see the references
at the beginning of this section.
One can see from the proof
that the flatness of $\fscon$
encodes many of the axioms
of Frobenius manifolds.

(\ref{eq:hbar-direction})
can be considered as
a family of meromorphic connections
on $\bP^1$
parametrized by $M$,
and the flatness of $\fscon$
means that
this family is isomonodromic.
One can see from (\ref{eq:hbar-direction})
that this meromorphic connection
has a regular singularity
at infinity
and an irregular singularity
of Poincar\'{e} rank 1
at the origin.

\begin{definition}
 A point $p$ on a Frobenius manifold $M$
 is semisimple
 if the product structure $\circ$
 on the tangent space $T_p M$ of $M$
 at $p$
 is semisimple
 (i.e., if there are no nilpotent elements).
\end{definition}

By Dubrovin \cite{Dubrovin_G2DTFT},
there exists a local coordinate
$\{ u_\canonicala \}_{\canonicala=1}^N$
satisfying
\begin{eqnarray}
 \partial_\canonicala \circ \partial_\canonicalb
  = \delta_{\canonicala \canonicalb} \partial_\canonicala,
 \label{eq:canonical_idempotent} \\
 E = \sum_{\canonicala=1}^N u_\canonicala \partial_\canonicala,
 \label{eq:canonical_euler}
\end{eqnarray}
where
$\partial_\canonicala
 = \frac{\partial}{\partial u_\canonicala}$,
in the neighborhood of a semisimple point.
Since
\begin{eqnarray*}
 g(\partial_\canonicala, \partial_\canonicala)
  &=& g(\partial_\canonicala, \partial_\canonicalb
           \circ \partial_\canonicalb) \\
  &=& g(\partial_\canonicala \circ \partial_\canonicalb,
           \partial_\canonicalb),
\end{eqnarray*}
$\partial_\canonicala$ and 
$\partial_\canonicalb$ are orthogonal
for $\canonicala \neq \canonicalb$.
This coordinate
$\{ u_\canonicala \}_{\canonicala=1}^N$
is called the {\em canonical coordinate},
which is unique up to permutations.
By
(\ref{eq:canonical_idempotent}) and
(\ref{eq:canonical_euler}),
$\{ u_\canonicala \}_{i=1}^N$ is characterized
as the set of eigenvalues of $\scU$.

Let
$$
 f_\canonicala = \Delta_\canonicala^{-1/2}
  \partial_\canonicala,
 \quad i=1, \ldots, N,
$$
be the normalized canonical vector field,
where
$$
 \Delta_\canonicala
  = g(\partial_\canonicala, \partial_\canonicala),
$$
and define an $N \times N$ matrix
$\Psi_{\flata \canonicala}$ by
$$
f_\canonicala
  = \sum_{\canonicala = 1}^N
     \Psi_{\flata \canonicala}  \frac{\partial}{\partial t_\flata}.
$$
$\Psi_{\flata \canonicala}$
is the coordinate transformation matrix
from the normalized canonical coordinate
to the flat coordinate.
%
The operator
$\scU$
defined in (\ref{eq:U})
is diagonal
in the normalized canonical coordinate:
$$
 \scU(\partial_\canonicala)
 = \sum_{\canonicalb=1}^N u_\canonicalb \partial_\canonicalb \circ
    \partial_\canonicala
 = u_\canonicala \partial_\canonicala.
$$
Let $U,V$ be the matrices
representing
the operators $\scU, \scV$
in the normalized canonical coordinate.
Then
$$
 U = \diag(u_1, \ldots, u_N)
$$
and (\ref{eq:hbar-direction}) becomes
\begin{equation} \label{eq:nc_hbar-direction}
 \partial_\hbar Y + \frac{1}{\hbar^2} U Y - \frac{1}{\hbar} V Y = 0.
\end{equation}

By Dubrovin
\cite{Dubrovin_PT2DTFT},
Lemma 4.3.,
there exists a unique
$N \times N$-matrix-valued formal series
$$
 R(\hbar) = 1 + R_1 \hbar + R_2 \hbar^2 + \cdots
$$
satisfying
$$
 R^{t}(\hbar) R(-\hbar) = 1
$$
such that
$$
 Y = R(\hbar) \widetilde{Y}
$$
transforms (\ref{eq:nc_hbar-direction})
to
\begin{equation} \label{eq:gt_hbar-direction}
 \hbar \partial_{\hbar} \widetilde{Y} + U \widetilde{Y} = 0.
\end{equation}
Here, $\bullet^t$ denotes the transpose
of a matrix.
Since (\ref{eq:gt_hbar-direction}) has a fundamental solution
of the form
$$
 \widetilde{Y} = \exp[U / \hbar],
$$
(\ref{eq:hbar-direction})
has a formal fundamental solution
of the form
\begin{equation} \label{eq:formal_solution}
 \Phi_\formal = \Psi R(\hbar) \exp[U / \hbar]
\end{equation}
in the flat coordinate.

\begin{definition}
For $0 \leq \phi < \pi$,
a straight line
$l = \{\hbar \in \bCx \suchthat \arg(\hbar) = \phi, \phi - \pi \}$
passing through the origin
is called {\em admissible}
if the line through $u_\canonicala$ and $u_\canonicalb$
is not orthogonal to $l$
for any $\canonicala \neq \canonicalb$.
\end{definition}

Fix such a line, and choose
a small enough number $\epsilon > 0$
so that any line passing through the origin
with angle between $\phi-\epsilon$ and $\phi+\epsilon$
is admissible.
\begin{figure}[h] \label{fg:admissible_line}
 \centering
 \psfrag{Dl}{$D_\lt$}
 \psfrag{Dr}{$D_\rt$}
 \psfrag{D-}{$D_-$}
 \psfrag{phi}{$\phi$}
 \psfrag{l}{$l$}
 \includegraphics{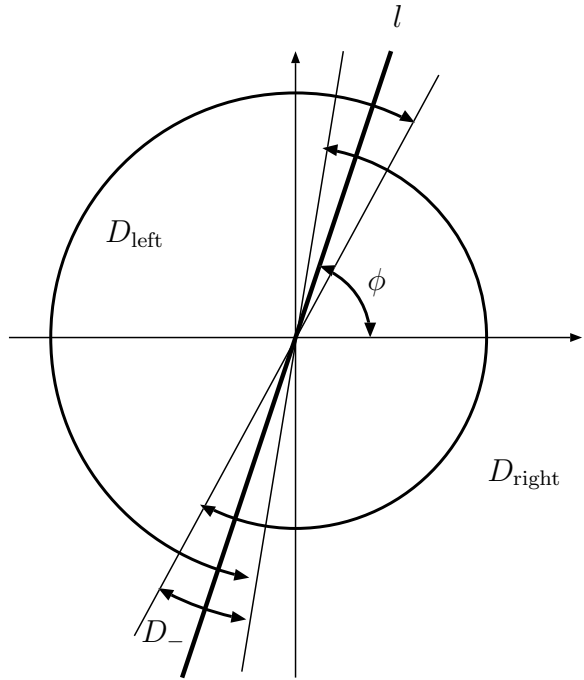}
 \caption{An admissible line and angular domains}
\end{figure}
Define
angular domains $D_{\rt}$, $D_{\lt}$ and $D_-$ by
\begin{eqnarray} 
D_{\rt} & = & \{\hbar \in \bCx \suchthat
 \phi - \pi - \epsilon < \arg(\hbar) <\phi + \epsilon \},
 \nonumber \\
D_{\lt} & = & \{\hbar \in \bCx \suchthat
 \phi - \epsilon < \arg(\hbar) <\phi +\pi + \epsilon \}, \\
D_- & = & \{\hbar \in \bCx \suchthat
 \phi - \pi - \epsilon < \arg(\hbar) < \phi - \pi + \epsilon \}.
 \nonumber
\end{eqnarray}
See Figure \ref{fg:admissible_line}.
Since the singularity at the origin
is irregular,
the formal solution $\Phi_\formal(\hbar)$
does not converge.
Nevertheless,
by \cite{Dubrovin_PT2DTFT} Theorem 4.2.,
there exist unique solutions
$\Phi_\rt(\hbar)$ and $\Phi_\lt(\hbar)$,
defined on the angular domains $D_\rt$ and $D_\lt$
respectively,
which asymptote
to the same formal solution:
$$
\Phi_{\rt/\lt} \sim \Phi_{\small \mathrm{formal}}
\ \ \mbox{as $\hbar \rightarrow 0$ in $D_{\rt/\lt}$}.
$$
Since these two solutions
satisfy the same linear differential equation
on $D_-$,
there exists a matrix $S$
independent of $\hbar$
such that
\begin{equation*}
 \Phi_\lt(\hbar) = \Phi_\rt(\hbar) S,
 \qquad 
 \quad \hbar \in D_-.
\end{equation*}
This matrix $S$ is called the Stokes matrix.

\section{Gromov-Witten invariants}
 \label{sc:Gromov-Witten_invariants}

Let $X$ be a smooth projective variety
over $\bC$.
For simplicity,
let us assume
that $H^k(X;\bC)=0$
for odd $k$.

\begin{definition}
$(f,C,p_1,\ldots,p_n)$
is a stable map
to $X$
of genus $g$
with $n$ marked points
if $C$ is a complete curve
of arithmetic genus $g$
with at worst ordinary double points,
$(p_1,\ldots,p_n)$
is an ordered set of $n$ points
on $C$,
$f:C \rightarrow X$
is a regular map
and the automorphism group
of $f$ is finite.
In this case,
$f_*[C] \in H_2(X,\bZ)$
is called the degree
of $f$.
\end{definition}

Let $\M_{g,n}(X;\degreed)$ be
the moduli space of
stable maps to $X$
of degree $\degreed$ and
genus $g$
with $n$ marked points.
We have the evaluation map
$$
\begin{array}{cccc}
 \ev_i: & \M_{g,n}(X;\degreed)
   & \longrightarrow & X \\
 & \rotatebox{90}{$\in$} & & \rotatebox{90}{$\in$} \\
 & (f,C,p_1,\ldots,p_n) & \longmapsto &
   f(p_i).
\end{array}
$$
Let $\scL_i$
be the line bundle on $\M_{g,n}(X;\degreed)$
whose fiber over
$(f,C,p_1,\ldots,p_n) \in \M_{g,n}(X;\degreed)$
is the cotangent space
$T_{p_i}^* C$
at the $i$-th marked point $p_i$.
The virtual dimension
of $\M_{g,n}(X;\degreed)$
is
$$
(1-g)(\dim_\bC X - 3) + \int_\degreed c_1(T X) + n,
$$
where
$c_1(T X)$ is the first Chern class
of the tangent bundle of $X$.
Given $n$ cohomology classes
$\gamma_1, \ldots, \gamma_n \in H^*(X;\bC)$
and non-negative integers
$d_1, \ldots, d_n$,
the Gromov-Witten invariant
is defined by
$$
 \langle \tau_{d_1} \gamma_1, \ldots,
  \tau_{d_n} \gamma_n \rangle_{g, \degreed}
  = \int_{[\M_{g,n}(X;\degreed)]^{\mathrm{virt}}}
     \prod_{i=1}^n
      \left(
       c_1(\scL_i)^{d_i} \cup \ev_i^*(\gamma_i)
      \right).
$$
Here,
$[\M_{g,n}(X;\degreed)]^{\mathrm{virt}}$
is the {\em virtual fundamental class}
\cite{Behrend-Fantechi},
\cite{Li-Tian_virtual}.

%

Let $\{ T_\flata \}_{\flata=0}^{N-1}$
be a basis of $H^*(X;\bC)$
and let $\{ t_\flata \}_{\flata=0}^{N-1}$
be the corresponding coorinate
of $H^*(X;\bC)$.
Set $\deg t_\flata = \deg T_\flata = k$
when $T_\flata \in H^k(X;\bC)$.
Assume that
the basis $\{ T_\flata \}_{\flata=0}^{N-1}$ are taken
so that $T_0$ is the basis of $H^0(X;\bC)$.
Let
$$
 \gamma = \sum_{\flata = 0}^{N-1} t_\flata T_\flata.
$$
The following formal power series
$\Phi(\gamma)$
in $\{t_i\}_{i=0}^{N-1}$
is called
the {\em Gromov-Witten potential}:
$$
 \Phi(\gamma) = \sum_{n=0}^\infty \sum_{\degreed \in H_2(X,\bZ)}
  \frac{1}{n!}
   \langle
    \underbrace{\gamma, \ldots, \gamma}_{n \ \text{times}}
   \rangle_{0, \degreed}. 
$$
Since
$$
 \frac{1}{n!}
   \langle
    \underbrace{\gamma, \ldots, \gamma}_{n \ \text{times}}
   \rangle_{0, \degreed}
 = \sum_{|\alpha| = n}
    \langle
     \underbrace{T_0, \ldots, T_0}_{\alpha_0 \ \text{times}},
     \ldots,
     \underbrace{T_{N-1}, \ldots, T_{N-1}}_{\alpha_{N-1} \ \text{times}}
    \rangle_{0, \degreed}
   \frac{t^\alpha}{\alpha !}
$$
where
$\alpha = (\alpha_0 , \ldots, \alpha_{N-1})$
is a multi-index,
$|\alpha|=\alpha_0+\cdots+\alpha_{N-1}$,
$\alpha!=\alpha_0! \cdots \alpha_{N-1}!$
and
$t^\alpha = t_0^{\alpha_0} \cdots t_{N-1}^{\alpha_{N-1}}$,
$\Phi$
is the generating function of
the Gromov-Witten invariants.

\begin{theorem}
When the Gromov-Witten potential $\Phi$ converges
in some domain $D$ in $H^*(X;\bC)$,
$D$ has a structure of a Frobenius manifold
with $\Phi$ as a potential function.
\end{theorem}

See the references at the beginning of section
\ref{sc:Frobenius_manifolds}
for details.
The above theorem means that
although the moduli space of stable maps
of genus zero
are divided into connected components
by the degree and the number of marked points,
the Gromov-Witten invariants
as a whole
have a strong structure
so that together
they turn
$H^*(X;\bC)$
into a Frobenius manifold.
A part of the Frobenius structure
is a product structure
on the tangent bundle.
Since
$T H^*(X;\bC) \cong H^*(X;\bC) \times H^*(X;\bC)$,
it is a product structure
on $H^*(X;\bC)$
parametrized by
$H^*(X;\bC)$ itself.
By the {\em Point Mapping Axiom}
in Gromov-Witten theory,
this product structure is
a deformation of the usual cup product
on the cohomology ring,
and
is called the {\em quantum cohomology ring}.
The metric $g$
is given by the Poincar\'{e} pairing
on the cohomology group.
The unit vector field is
$$
 e = \frac{\partial}{\partial t_0},
$$
the Euler vector field is
\begin{equation} \label{eq:gw-euler}
 E = \sum_{\flata=0}^{N-1}
      (1-\frac{\deg t_\flata}{2}) t_\flata
        \frac{\partial}{\partial t_\flata}
     + \sum_{\deg t_\flatb = 2}
        r_\flatb \frac{\partial}{\partial t_\flatb}, 
\end{equation}
and the charge is
$$
 D=2-\dim_\bC X.
$$
Here, $r_\flatb \in \bZ$ is defined by
\begin{equation} \label{eq:rflatb}
 c_1(T X) = \sum_{\deg t_\flatb = 2} r_\flatb T_\flatb.
\end{equation}

%

The WDVV equation
(\ref{eq:WDVV})
is a highly non-trivial differential equation
for the Gromov-Witten potential
$\Fpotential$,
which allows us
to compute
all the genus-zero Gromov-Witten invariants
from a few initial conditions
in special cases.
Fro the Gromov-Witten invariants
of rational surfaces,
see, e.g.,
\cite{Caporaso-Harris},
\cite{Crauder-Miranda},
\cite{Francesco-Itzykson},
\cite{Gathmann_GWIBU},
\cite{Goettsche-Pandharipande},
\cite{Kontsevich-Manin_reconstruction}
and references therein.

%

\section{Mirror symmetry for toric manifolds}
\label{sc:toric_mirror}

In this section,
we review the mirror symmetry for toric varieties.
See also \cite{Batyrev_QCRTM},
\cite{Givental_EGWI},
\cite{Givental_MTTCI},
\cite{Iritani_QDGMT}.

Let $\Sigma$ be a complete fan
in $N \cong \bZ^\dimfan$
and $X_\Sigma$ be the corresponding toric variety.
We denote the primitive generators
of one-dimensional cones of $\Sigma$
by $\{v_\edgea\}_{\edgea=1}^\noedges$,
$
 v_\edgea = (v_{\edgea 1},\ldots,v_{\edgea \dimfan})
   \in \bZ^\dimfan.$
Here,
$\picn$
is the Picard number
of $\XSigma$.
If we put $M=\Hom(N, \bZ)$,
we have an exact sequence
\begin{equation*}
 0 \rightarrow M \rightarrow \bZ^\noedges \rightarrow
    H^2(\XSigma;\bZ) \rightarrow 0,
     \label{eq:htwo} 
\end{equation*}
where
the homomorphism
$M \rightarrow \bZ^\noedges$
is defined by
$$
M \ni m \mapsto (m(v_\edgea))_{\edgea=1}^\noedges \in \bZ^\noedges.
$$
Choose a basis
$\{ T_\flata \}_{\flata=1}^\picn$
of $H^2(X;\bZ)$
and let $t = (t_1, \ldots, t_\picn)$
be the corresponding coordinate
of $H^2(X;\bZ)$.
One-dimensional cones in
$\Sigma$
correspond to line bundles
on $\XSigma$,
and let $\toricu_\edgea$ be
the first Chern class
of the line bundle
corresponding to
the one-dimensional cone
generated by
$v_\edgea$.
Define an $(\noedges) \times \picn$-matrix
$m_{\edgea \flata}$ by
\begin{equation} \label{eq:miflata}
 \toricu_\edgea = \sum_{\flata=1}^\picn m_{\edgea \flata} T_\flata.
\end{equation}

Now,
following Givental,
introduce
an $H^*(\XSigma;\bC)$-valued
formal function
$I(\motonot_0, \motonot; \hbar)$
of
$\motonot_0 \in H^0(\XSigma; \bC)$,
$\motonot \in H^2(\XSigma;\bC)$ and
$\hbar \in \bCx$ by
\begin{equation} \label{eq:defI}
 I(\motonot_0, \motonot; \hbar)
 = e^{(\motonot_0 + T \motonot) / \hbar}
 \sum_{\degreed \in H_2(\XSigma; \bZ)_{\mathrm{eff}}}
  q^\degreed
   \prod_{\edgea=1}^\noedges
    \frac{\prod_{\Irindex=-\infty}^0 (\toricu_\edgea + \Irindex \hbar)}
     {\prod_{\Irindex=-\infty}^{\langle \degreed, \toricu_\edgea \rangle}
       (\toricu_\edgea + \Irindex \hbar)}. 
\end{equation}
Here,
$T \motonot = T_1 \motonot_1 + \cdots + T_\picn \motonot_\picn$,
$H_2(\XSigma; \bZ)_{\mathrm{eff}} \subset H_2(\XSigma; \bZ)$
is the cone
spanned by the homology class of curves
and
$$
 q^\degreed = q_1^{\pair{\degreed}{T_1}}
               \cdots q_\picn^{\pair{\degreed}{T_\picn}},
$$
where
$q_\flata=\exp(t_\flata)$,
$\flata=1,\ldots,\picn$ and
$\langle \bullet, \bullet \rangle$
is the pairing
between homology and cohomology.

It follows from (\ref{eq:defI})
that the $I$-function
satisfies the following
set of hypergeometric equations
of the GKZ (Gelfand-Kapranov-Zelevinski) type:

\begin{equation} \label{eq:GKZ}
 D_\degreed I = 0, \quad
  \degreed \in H_2(\XSigma;\bZ),
\end{equation}
where
$$
 D_\degreed = q^\degreed
  \prod_{\pair{\toricu_\edgea}{\degreed} < 0}
   \prod_{\Irindex=0}^{- \pair{\toricu_\edgea}{\degreed} - 1}
    U_{\edgea \Irindex}
  -
  \prod_{\pair{\toricu_\edgea}{\degreed} > 0}
   \prod_{\Irindex=0}^{\pair{\toricu_\edgea}{\degreed} - 1}
    U_{\edgea \Irindex},
$$
$$
 U_{\edgea \Irindex}
  = \sum_{\flata=1}^\picn m_{\edgea \flata}
     \hbar \frac{\partial}{\partial t_\flata} - \Irindex \hbar.
$$
From now on, we assume that
the anticanonical bundle $-K_{\XSigma}$
of the toric variety $\XSigma$ is nef.
Note that
although $I$-function is an infinite series
in both positive and negative powers of $\hbar$
for general $\XSigma$,
it contains only non-positive powers of $\hbar$
when $-K_{\XSigma}$ is nef.
This follows from the fact that
$
\exp[-(t_0 + T t)/\hbar] I
$
is homogeneous of degree zero
with respect to the degree
assigned as
\begin{equation} \label{eq:degree_assignment}
 \deg \toricu_\canonicala = \deg \hbar = 1, \ 
 \deg( q^\degreed ) = \langle \degreed, c_1(T \XSigma) \rangle.
\end{equation}
Therefore,
the coefficient of a positive power of $\hbar$
must have a negative degree,
which is impossible
when $-K_{\XSigma}$ is nef.

I learned the following Lemma \ref{lemma:nef_convergence}
from H. Iritani:

\begin{lemma} \label{lemma:nef_convergence}
When the anticanonical bundle $-K_{\XSigma}$
of a toric variety $\XSigma$ is nef,
there exists a domain $D \subset H^2(\XSigma; \bC)$
such that
$I(\motonot_0, \motonot; \hbar)$ converges
for any $\motonot_0 \in H^0(\XSigma; \bC)$,
$\motonot \in D$
and $|\hbar|>1$
\end{lemma}
\proof
Fix a norm $\| \bullet \|$
on $H^*(\XSigma;\bC)$
such that
$\| \omega \tau \| \leq \| \omega \| \| \tau \|$.
Let $C_1$ be a positive number
satisfying
$$
 \sup_{|\hbar|=1, \ \Irindex \ge 1, \ \edgea=1, \ldots, \noedges}
   \left(
    \left\| \toricu_\edgea \right\|,
    \left\| 1+\frac{\toricu_\edgea}{\Irindex \hbar} \right\|,
    \left\| \left( 1 + \frac{\toricu_\edgea}{\Irindex\hbar}
            \right)^{-1} \right\|
   \right)
  \le C_1.
$$
Then
for $\hbar \in \bC$ such that $|\hbar| = 1$,
\begin{align*}
 \left\| \prod_{\edgea=1}^\noedges
  \frac{\prod_{\Irindex=-\infty}^0 (\toricu_\edgea + \Irindex \hbar)}
    {\prod_{\Irindex=-\infty}^{\langle \degreed, \toricu_\edgea \rangle}
     (\toricu_\edgea + \Irindex \hbar)}
 \right\|
& \le
\frac{\prod_{\langle \degreed, \toricu_\edgea \rangle <0}
             |\langle \degreed, \toricu_\edgea \rangle|!}
     {\prod_{\langle \degreed, \toricu_\edgea \rangle >0}
             \langle \degreed, \toricu_\edgea \rangle!}
 C_1^{\sum_{\edgea=1}^\noedges |\langle \degreed, \toricu_\edgea \rangle|}  \\
& \le \frac{ C_2^{|\degreed|} }
           { \langle \degreed, c_1(X) \rangle !}
\end{align*}
Here,
$C_2$ is some positive number and
$|\degreed| = |\langle \degreed, T_1 \rangle|
   + \cdots + |\langle \degreed, T_\picn \rangle|$.
The second inequality follows from
\begin{align*}
 n_1!n_2!\cdots n_{k}! & \le (n_1+\cdots+n_k)!, \\
 \frac{1}{n_1! n_2! \cdots n_k!}
  &= \frac{1}{(n_1+\cdots+n_k)!}
      \binom{n_1+\cdots +n_k}{n_1,\cdots,n_k} \\
  & \le \frac{1}{(n_1+\cdots+n_k)!} k^{n_1+\cdots+n_k}, \\
 \frac{n!}{m!}
  & \le \frac{1}{(m-n)!} \text{\quad if } n \le m.
\end{align*}
and
$$
 \sum_{\edgea=1}^{\noedges} \langle \degreed, \toricu_\edgea \rangle
  = \langle \degreed, -K_\XSigma \rangle \geq 0
$$
since $\degreed$ is an effective class
and $-K_{\XSigma}$ is nef.
Therefore,
if we define
$
 D = \{ t \in H^2(\XSigma;\bC) \suchthat 
        q_\flata < 1 / C_2, \ 
         \flata = 1, \ldots, \picn. \},
$
the $I$-function converges
when $|\hbar|=1$
and $t \in D$.
Since the $I$-function contains
only non-positive powers of $\hbar$,
it also converges
when $|\hbar| > 1$
and $t \in D$.
\qed

Now, let us
introduce
another $H^*(\XSigma;\bC)$-valued
formal function
$J(t_0, t; \hbar)$
of $t_0 \in H^0(\XSigma;\bC)$,
$t \in H^2(\XSigma;\bC)$
and $\hbar \in \bCx$
as follows:
First,
define a set of
$H^*(\XSigma;\bC)$-valued formal functions
$\{ s_\flata \}_{\flata=0}^{N-1}$
of $\{t_\flata \}_{\flata=0}^{\picn}$
by
$$
 s_\flata = T_\flata + \sum_{n=0}^\infty \hbar^{-n-1}
  \sum_{\flatb=0}^{N-1}
   \sum_{k=0}^\infty
    \sum_{\degreed \in H_2(\XSigma; \bC)}
     \frac{1}{k!}
      \langle
       \tau_n T_\flata, T_\flatb,
       \underbrace{\gamma, \ldots, \gamma}_{k \ \text{times}}
      \rangle_{0,\degreed}
     T^\flatb,
$$
where
$$
 \gamma = \sum_{\flata = 0}^\picn t_\flata T_\flata.
$$
Here, $T^\flata$ is the Poincar\'{e} dual
of $T_\flata$:
$$
 g(T_\flata, T^\flatb) = \delta_{\flata \flatb}.
$$
Then it follows from the
{\em topological recursion relation}
that
\begin{equation} \label{eq:trr_flat}
 \hbar \frac{\partial s_\flata}{\partial t_\flatb}
  = T_\flatb \circ s_\flata,
 \quad \flatb = 0, 1, \ldots, \picn.
\end{equation}
See, e.g., \cite{Cox-Katz} Proposition 10.2.1.
Note that
(\ref{eq:trr_flat})
is the $t_\flatb$-part of
the differential equation
for flat sections
of the first structure connection
(\ref{eq:X-direction}).
Now define
the {\em Givental's J-function} by
$$
 J = \sum_{\flata=0}^{N-1} g(s_\flata, 1) T^\flata.
$$

Since $H^*(\XSigma;\bC)$
is generated by
$H^2(\XSigma;\bC)$
as a ring,
there exists a set of polynomials
$\{ P_\flata(x_1, \ldots, x_\picn) \}_{\flata=0}^{N-1}$
of $\picn$ variables
such that
$$
 T_\flata = P_\flata(T_1, \ldots, T_\picn),
  \quad \flata = 0, 1, \ldots, N-1.
$$
Since
$$
 T_\flata \circ T_\flatb = T_\flata \cup T_\flatb + O(q)
$$
by the Point Mapping Axiom,
there exists a set
$\{ \widetilde{P}_\flata(x_1, \ldots, x_\picn)
 \}_{\flata = 0}^{N-1}$
of polynomials
whose coefficient is a formal power series in $q$
such that
$$
 \widetilde{P}_\flata(
      \hbar \frac{\partial}{\partial t_1} + T_1 \circ,
      \ldots,
      \hbar \frac{\partial}{\partial t_\picn} + T_\picn \circ) 1
  = T_\flata, \quad \flata = 0, 1, \ldots, N-1.
$$
Then, since
\begin{align*}
 P(\hbar \partial_1, \ldots, \hbar \partial_\picn) J
  &= \sum_{\flata=0}^{N-1}
       g \left(s_\flata,
         P(\hbar \frac{\partial}{\partial t_1} + T_1 \circ,
           \ldots,
           \hbar \frac{\partial}{\partial t_\picn} + T_\picn \circ
          ) 1
           \right) T^\flata \\
\end{align*}
for a polynomial $P(x_1, \ldots, x_\picn)$,
$$
 s_\flata
  = \sum_{\flatb = 0}^{N-1} g \left(
     \widetilde{P}_\flatb(\hbar \partial_1, \ldots, \hbar \partial_\picn) J,
     T_\flata \right) T^\flatb.
$$

To sum up,
we have introduced
two $H^*(\XSigma;\bC)$-valued
formal functions
$I(t_0, t; \hbar)$ and $J(t_0, t; \hbar)$.
The $I$-function is defined
in terms of the combinatorial data
of the fan $\Sigma$ defining $\XSigma$
and satisfies a set of hypergeometric
differential equations
(\ref{eq:GKZ}) of GKZ type.
The $J$-function
is defined in terms of the Gromov-Witten invariants
of $\XSigma$
and encodes the information
of the solution
to the equation
(\ref{eq:trr_flat}),
hence the structure of the product $\circ$.

The mirror theorem by Givental
(see \cite{Givental_MTTCI} Theorem 0.2,
Proposition 6.4, and the discussion
at the end of section 7)
gives the relation between
the $I$-function
and the $J$-function
as follows:
Since $-K_\XSigma$ is nef,
all the variables $\motonoq_\flata$'s
have non-negative degrees,
see (\ref{eq:degree_assignment}).
It follows from the homogeneity
of
$\exp[- (\motonot_0 + T \motonot) / \hbar]
  I(\motonot_0, \motonot; \hbar)$
that the $I$-function has the form
\begin{equation} \label{eq:expansion_I}
 I(\motonot_0, \motonot;\hbar)
  = e^{(\motonot_0 + T \motonot)/\hbar}
     \left( \Izero(\motonoq)
        + (\Ionezero(\motonoq) + \Ioneone(\motonoq)) \hbar^{-1}
      + O(\hbar^{-2}) \right),
\end{equation}
where
$\Izero(\motonoq)$ and $\Ionezero(q)$
are $H^0(\XSigma; \bC)$-valued
and
$\Ioneone(q)$ is $H^2(\XSigma; \bC)$-valued.
Then, the $J$-function is given by
\begin{align*}
 J(\atonot_0, \atonot; \hbar)
  &= I(\motonot_0, \motonot; \hbar) / \Izero(\motonoq) \\
  &= \exp[(\motonot_0 + T \motonot)/\hbar]
      \left(
       1 + \Izero(\motonoq)^{-1}
        \left(
         \Ionezero(\motonoq)
          + \Ioneone(\motonoq)
        \right) \hbar^{-1}
       + O(\hbar^{-2})
      \right) \\
  &= \exp
      \left[
       \left(
        \left(
         \motonot_0 + \Izero(\motonoq)^{-1} \Ionezero(e^\motonoq)
        \right)
        + \left(
           \motonot + \Izero(\motonoq)^{-1} \Ioneone(\motonoq)
          \right)
       \right) /\hbar
      \right]
      (1 + O(\hbar^{-2})) \\
  &= \exp \left[
           \left( \atonot_0 + T \atonot \right) /\hbar
          \right]
       (1 + O(\hbar^{-2})).
\end{align*}
Here,
we have performed a
homogeneous coordinate transformation
$$
 \atonot_0
  = \motonot_0 + \Izero(\motonoq)^{-1} \Ionezero(\motonoq),
 \quad
 \atonot = \motonot + \Izero(\motonoq)^{-1} \Ioneone(\motonoq),
$$
called the {\em mirror transformation}.
Since the $I$-function is convergent,
the mirror transformation
and the $J$-function are also convergent.

\section{Stationary-phase integral and Stokes matrix}
 \label{sc:spint}
Let $v_0 = (0, \ldots, 0) \in \bZ^\dimfan$,
and $\{ v_\edgea \}_{\edgea=1}^\noedges \subset \bZ^\dimfan$
be the set of primitive generators
of one-dimensional cones of
the fan $\Sigma$ in $\bZ^\dimfan$
defining a toric variety $\XSigma$
as in section \ref{sc:toric_mirror}.
We assume that
the anticanonical bundle
$-K_{\XSigma}$ of $\XSigma$ is nef.
Define a Laurent polynomial
$W$ by
\begin{equation} \label{eq:toric_mirror} 
 W(x_1,\ldots,x_n)
  = \sum_{\edgea=0}^\noedges \wa_\edgea
     x_1^{v_{\edgea 1}} \cdots x_r^{v_{\edgea r}}.
\end{equation}
$W$ defines a regular function
on the algebraic torus
$(\bCx)^\dimfan$.
By Kouchnirenko \cite{Kouchnirenko},
the number $N$
of critical points
of $W$ is
$\dimfan !$ times the volume
of the Newton polygon
of $W$
(i.e., the convex hull of $\{v_\edgea \}_{\edgea =1}^\noedges$)
when $\wa_\edgea \in \bC$, $\edgea=0,\ldots,\noedges$,
are general enough.
This $N$ is equal to
the rank 
of the total cohomology ring
$H^*(\XSigma,\bC)$
of $\XSigma$.
Let $\{p_\canonicala\}_{\canonicala=1}^N$
be the set of critical points of $W$.
They are functions of
$\{ \wa_\edgea \}_{\edgea = 0}^\noedges$.
Fix a complete K\"{a}hler metric
on the algebraic torus $(\bCx)^\dimfan$.
For a general $\hbar$,
let $\Gamma_\canonicala(\hbar)$
be the descending Morse cycle for
$\Re(W/\hbar)$.
The image of $\Gamma_\canonicala(\hbar)$
by $W$
is a half-line
starting from the critical value
$W(p_\canonicala)$,
and the fiber
above a point $p$
on this half-line
is the cycle
in $W^{-1}(p)$
which vanishes
at $p_\canonicala$
by the parallel transport
along this half-line.
Following Givental
\cite{Givental_EGWIGMC},
for $\canonicala=1,\ldots,N$,
consider the integral
\begin{equation} \label{eq:spint}
 \spint_\canonicala
  = (\pi \hbar)^{- \dimfan / 2}
    \int_{\Gamma_\canonicala(\hbar)} e^{\frac{1}{\hbar} W(x)}
     \frac{d x_1 \wedge \cdots \wedge d x_n}{x_1 \cdots x_\dimfan}.
\end{equation}
This integral is invariant
under the natural action
of the torus $(\bCx)^\dimfan$
on itself defined by
$$
 (\bCx)^\dimfan \ni (\alpha_1,\ldots,\alpha_\dimfan) :
  (x_1,\ldots,x_\dimfan) \mapsto
     (\alpha_1 x_1, \ldots, \alpha_r x_\dimfan).
$$
Therefore,
if we define the variables
$\{ t_\flata \}_{\flata=0}^\picn$
by $t_0 = \wa_0$ and
$$
 \exp({t_\flata})
  = \wa_1^{m_{1 \flata}} \cdots \wa_\noedges^{m_{\noedges, \flata}},
 \quad \flata = 1, \ldots, \picn,
$$
where $m_{\edgea \flata}$'s are defined
in (\ref{eq:miflata}),
$\spint_\canonicala$'s
depend on $\{ \wa_\edgea \}_{\edgea=0}^\noedges$
only through $\{ t_\flata \}_{\flata=0}^\picn$.


It is easy to see
that $\spint_\canonicala$'s
satisfy the same differential equations
(\ref{eq:GKZ})
as the $I$-function.
Since the solution space to
these differential equation
is an $N$-dimensional $\bC$-vector space
and $\{ \spint_\canonicala \}_{\canonicala=1}^N$ are
linearly independent,
the set of components of the $I$-function and
$\{ \spint_\canonicala \}_{\canonicala=1}^N$ differ
only by an $\hbar$-dependent
linear transformation.

The $\hbar$-dependence
of $\spint_\canonicala$ is simple.
If we assign
degree $1$ to variables $t_0$, $\hbar$
and $\{ x_i \}_{i=1}^\dimfan$,
and degree $r_\flata$
to the variable $\exp(t_\flata)$
for $\flata=1, \ldots, \picn$
where
$r_\flata$'s are defined
in (\ref{eq:rflatb}),
$\spint_\canonicala$ is homogeneous
of degree $-\dimfan/2$:
\begin{equation} \label{eq:homogeneity}
 \left[ \hbar \frac{\partial}{\partial \hbar}
  + t_0 \frac{\partial}{\partial t_0}
  + \sum_{\flata=1}^\picn
     r_\flata \frac{\partial}{\partial t_\flata}
  + \frac{\dimfan}{2} \right] \spint_\canonicala
  = 0. 
\end{equation}

The important point is that
the $\hbar$-part
(\ref{eq:hbar-direction})
of the differential equation
for the $H^{2 \dimfan}(\XSigma;\bC)$-component
of a flat section
of the first structure connection
is precisely the homogeneity condition
(\ref{eq:homogeneity}):
Assume that a section $Y$ of $\pi^* \scT_M$
satisfies
$$
 \fscon_X Y = \nabla_X Y
  - \frac{1}{\hbar} X \circ Y
  = 0
$$
for any section $X$ of $\pi^* T_M$.
Then
\begin{align*}
 \hbar \fscon_{\partial_\hbar} Y
  &= \hbar \partial_\hbar Y
      + \frac{1}{\hbar} \scU(Y)
      - \scV(Y) \\
  &= \hbar \partial_\hbar Y
      + \frac{1}{\hbar} E \circ Y
      - \scV(Y) \\
  &= \hbar \partial_\hbar Y
      + \nabla_E Y
      - \scV(Y).
\end{align*}
Note that
for the flat coordinates
$\{ t_\flata \}_{\flata=0}^{N-1}$,
\begin{align*}
 \scV(\partial_\flata)
  &= \nabla_{\partial_\flata} E - \frac{D}{2} \partial_\flata \\
  &= \nabla_{\partial_\flata} 
      \left(
        \sum_{\flatb=0}^{N-1}
           (1-\frac{\deg t_\flatb}{2}) t_\flatb
            \partial_\flatb
        + \sum_{\deg t_\flatc = 2}
            r_\flatc \partial_\flatc
      \right)
      - \frac{2 - \dimfan }{2} \partial_\flata \\
  &= \left( \frac{\dimfan - \deg t_\flata}{2} \right)
      \partial_\flata.
\end{align*}
Since $\nabla_E$ counts the degree,
this differential equation
requires
the $H^k(\XSigma; \bC)$-component
of the flat section
of the first structure connection
to be homogeneous of degree
$(\dimfan - k) / 2$.

Now, recall that by Givental's mirror theorem,
$\Izero(\motonoq)^{-1} I(\motonot_0, \motonot; \hbar)$
coincides with the $J$-function,
which have the same dependence
on $\{ t_\flata \}_{\flata=0}^\picn$
with the $H^{2 \dimfan}(\XSigma;\bC)$-components
of the flat section
for the first structure connection.
Then, since the stationary-phase integrals
(\ref{eq:spint})
have the same dependence
on $\{ t_\flata \}_{\flata=0}^\picn$
as the $I$-function
and satisfies the desired homogeneity condition
(\ref{eq:homogeneity}),
$\Izero(\motonoq)^{-1} \spint_\canonicala(\motonot_0, \motonot; \hbar)$
gives the $H^{2 \dimfan}(\XSigma;\bC)$-component
of a flat section
for the first structure connection
for $\canonicala = 1, \ldots, N$.

By the saddle-point approximation,
\begin{equation} \label{eq:saddle-point_method}
 \Izero(\motonoq)^{-1} \spint_\canonicala
  \sim
  \frac{\Izero(\motonoq)^{-1}}
       {p_{\canonicala,1}\cdots p_{\canonicala,n}}
  \left[ \det \left(
   \frac{\partial^2 W}{\partial x_l \partial x_m} (p_\canonicala)
              \right)_{l,m} \right]^{-\frac{1}{2}}
  \left(
   1+O(\hbar)
  \right)
  e^{\frac{1}{\hbar}W(p_\canonicala)}
\end{equation}
as $|\hbar| \rightarrow 0$
with $\arg \hbar$ fixed,
where
$
 p_\canonicala = (p_{\canonicala,1}, \cdots,
   p_{\canonicala,n}) \in (\bCx)^\dimfan
$
is the $\canonicala$-th critical point
of $W$.
This gives the asymptotic expansion
appearing in (\ref{eq:formal_solution}).

Since the Stokes matrix acts
on all the components of
the fundamental solution
equally,
it is enough to compute
its action on the $H^{2 \dimfan}(\XSigma;\bC)$-component
computed above.
The integration cycle
$\Gamma_{\canonicala}(\hbar)$
in (\ref{eq:spint})
undergoes a discontinuous change
when $\hbar$ crosses the line
such that
the half-line
starting from the critical value
$u_\canonicala$
in the direction of $- \hbar^{-1}$
passes through another critical value
$u_\canonicalb$.
Therefore,
$\spint_\canonicala$
is not holomorphic
on such a line.
Since $\spint_\canonicala$
has a monodromy,
it is impossible to obtain
a holomorphic function
on the whole $\hbar$-plane
by analytic continuation,
but one can obtain a holomorphic function
on the two angular domains
$D_\lt$ and $D_\rt$
defined in Section \ref{sc:Frobenius_manifolds}
as follows:
Consider the local system
on $\bCx$
whose fiber over $\hbar \in \bCx$ is
the relative homology group
$H_n \left( (\bC^\times)^n, \Re(W/\hbar) \ll 0; \bZ \right)$.
Define
two sets of sections
$\{ \Gamma_{\canonicala, \lt}(\hbar) \}_{\canonicala=1}^N$
and
$\{ \Gamma_{\canonicala, \rt}(\hbar) \}_{\canonicala=1}^N$
of this local system
on $D_\lt$ and $D_\rt$
by the condition that
they coincide with the classes
defined by
the descending Morse cycles
at $\hbar = \exp[ (\phi + \pi / 2) \sqrt{-1}]$
and $\hbar = \exp[ (\phi - \pi / 2) \sqrt{-1}]$
respectively;
\begin{align*}
 \Gamma_{\canonicala, \lt}(\exp[ (\phi + \pi / 2) \sqrt{-1} ])
  &= \Gamma_{\canonicala}(\exp[ (\phi + \pi / 2) \sqrt{-1} ]), \\
  \Gamma_{\canonicala, \rt}(\exp[ (\phi - \pi / 2) \sqrt{-1} ])
  &= \Gamma_{\canonicala}(\exp[ (\phi - \pi / 2) \sqrt{-1} ]).
\end{align*}
Since $D_\lt$ and $D_\rt$ are simply-connected,
these conditions determine
$\Gamma_{\canonicala, \lt}(\hbar)$
and
$\Gamma_{\canonicala, \rt}(\hbar)$
uniquely.
Define $\spint_{\canonicala, \lt}$
(resp. $\spint_{\canonicala, \rt}$)
by
the integral (\ref{eq:spint})
with $\Gamma_{\canonicala, \lt}(\hbar)$
(resp. $\Gamma_{\canonicala, \rt}(\hbar)$)
as the integration cycle.
Then $\spint_{\canonicala, \lt}$
(resp. $\spint_{\canonicala, \rt}$)
are holomorphic
and
have the asymptotic behaviors
(\ref{eq:saddle-point_method})
as $\hbar \rightarrow 0$
in $D_\lt$
(resp. $D_\rt$).

\begin{figure}[htbp]
\begin{minipage}{0.5 \textwidth}
\centering
\psfrag{p1}{$W(p_1)$}
\psfrag{p2}{$W(p_2)$}
\psfrag{p3}{$W(p_3)$}
\psfrag{pN}{$W(p_N)$}
\psfrag{G1}{$\Gamma_1$}
\psfrag{G2}{$\Gamma_2$}
\psfrag{G3}{$\Gamma_3$}
\psfrag{GN}{$\Gamma_N$}
\psfrag{Re}{$\Re$}
\psfrag{Im}{$\Im$}
\includegraphics[width=\textwidth]{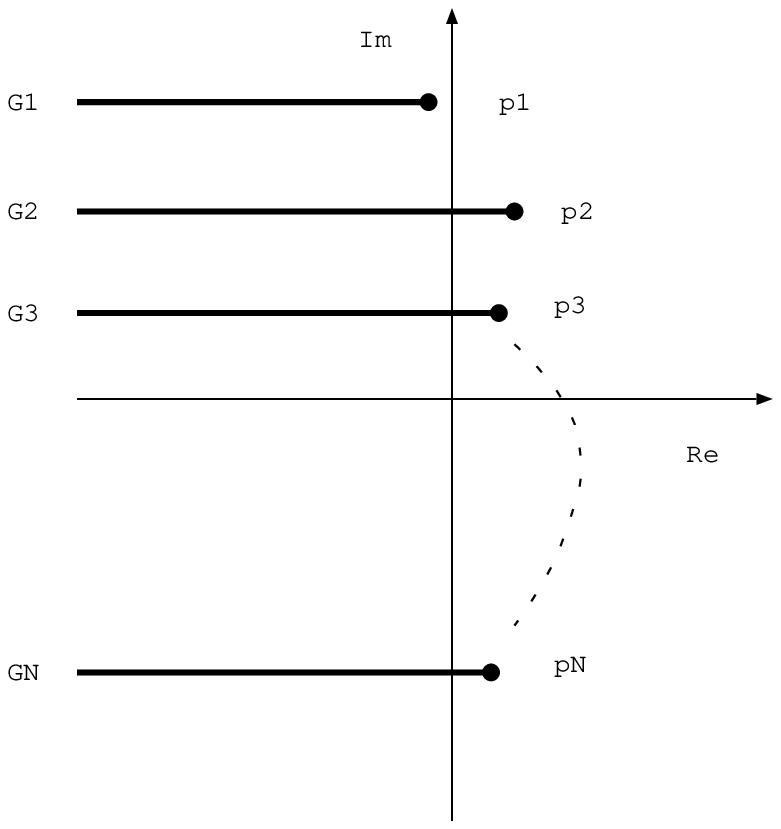}
\caption{
 Integration cycles
 at $\hbar = 1$
}
\label{fg:stokes1} 
\end{minipage}
\begin{minipage}{0.5 \textwidth}
\centering
\psfrag{p1}{$W(p_1)$}
\psfrag{p2}{$W(p_2)$}
\psfrag{p3}{$W(p_3)$}
\psfrag{pN}{$W(p_N)$}
\psfrag{Gr1}{$\Gamma_{1, \rtt}$}
\psfrag{Gr2}{$\Gamma_{2, \rtt}$}
\psfrag{Gr3}{$\Gamma_{3, \rtt}$}
\psfrag{GrN}{$\Gamma_{N, \rtt}$}
\psfrag{Gl1}{$\Gamma_{1, \ltt}$}
\psfrag{Gl2}{$\Gamma_{2, \ltt}$}
\psfrag{Gl3}{$\Gamma_{3, \ltt}$}
\psfrag{GlN}{$\Gamma_{N, \ltt}$}
\psfrag{Re}{$\Re$}
\psfrag{Im}{$\Im$}
\includegraphics[width=\textwidth]{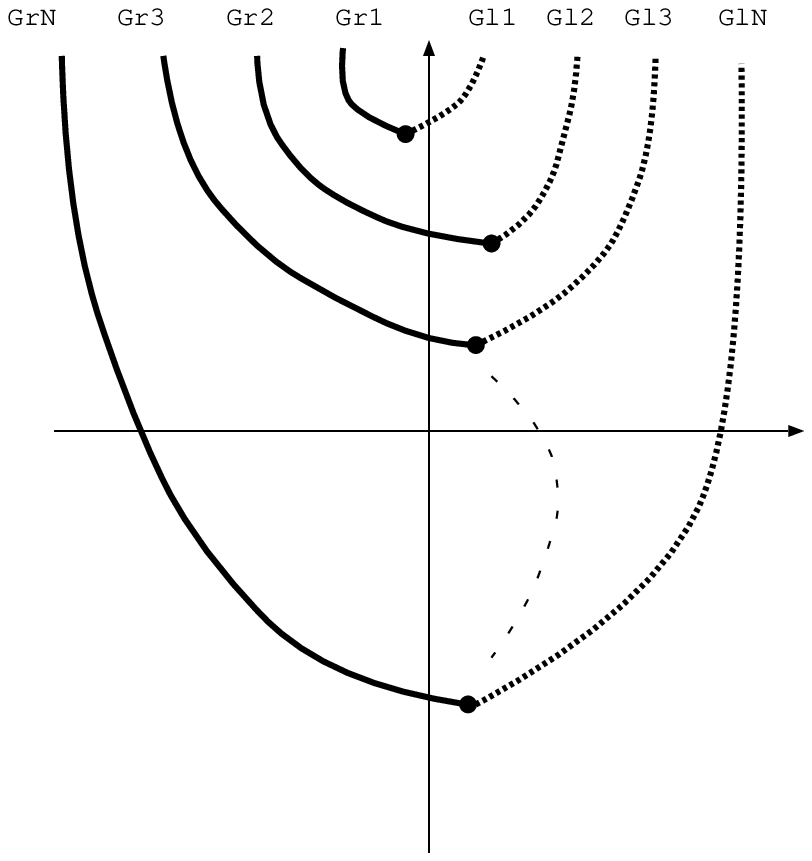}
\caption{
 Integration cycles
 at $\hbar = -\sqrt{-1}$
}
\label{fg:stokes2}
\end{minipage}
\end{figure}

Figure \ref{fg:stokes1}
shows the images $\Gamma_\canonicala$'s
of $\Gamma_{\canonicala, \rt}(\hbar)$'s
by $W$
at $\hbar = 1$.
Here,
we have assumed that
the imaginary axis
$\Re(\hbar) = 0$ is admissible,
and have chosen it
as an admissible line $l$.
Figure \ref{fg:stokes2}
shows
the images
$\Gamma_{\canonicala, \ltt}$'s
and $\Gamma_{\canonicala, \rtt}$'s
of
$\Gamma_{\canonicala, \lt}(\hbar)$'s
and $\Gamma_{\canonicala, \rt}(\hbar)$'s
by $W$
at $\hbar = - \sqrt{-1}$.

Since the integrand is single-valued,
the Stokes matrix is given
by the transformation matrix
between these two sets of cycles.
To describe it,
let us first choose
a regular value
$p \in \bC$
such that $\Re(p)$ is small enough
so that the line segment
$c_\canonicala$
from the $\canonicala$-th critical value
$u_\canonicala$
to $p$
is above $c_\canonicalb$
for any $\canonicalb > \canonicala$
(i.e., if $x \in c_\canonicala$,
$y \in c_\canonicalb$
and $\Re(x) = \Re(y)$,
then $\Im(x) > \Im(y)$).
Then,
the ordered set of cycles
$(C_\canonicala)_{\canonicala = 1}^N$
in $W^{-1}(p)$,
which vanish along the paths
$(c_\canonicala)_{\canonicala = 1}^N$
respectively,
forms a
{\em distinguished basis of vanishing cycles},
see e.g. \cite{Arnold--Gusein-Zade--Varchenko}.
It is obvious that
$$
\Gamma_{1, \lt}(\hbar)
 = \Gamma_{1, \rt}(\hbar). \\
$$
in $H_n((\bC^\times)^n, \Re[W/\hbar] \ll 0; \bZ)$
in the neighborhood of $\hbar = - \sqrt{-1}$.
As for $\Gamma_2(\hbar)$,
it follows from the Picard-Lefschetz formula
that
$$
 \Gamma_{2, \lt}(\hbar)
  = \Gamma_{2,\rt}(\hbar)
     - (C_1,C_2) \Gamma_{1, \rt}(\hbar),
$$
where
$(\bullet, \bullet)$
is the intersection form
in $H_\dimfan(W^{-1}(p); \bZ)$.
In the same way,
we have
\begin{eqnarray}
\Gamma_{\canonicala,\lt}(\hbar) 
 &=& \Gamma_{\canonicala,\rt}(\hbar)
   - (C_{\canonicala-1}, C_{\canonicala}) \Gamma_{\canonicala-1,\rt}(\hbar)
   - \cdots - (C_1, C_\canonicala) \Gamma_{1,\rt}(\hbar),
  \nonumber
\end{eqnarray}
by successive use of
the Picard-Lefschetz formula,
from which
\begin{eqnarray}
 S_{\canonicala \canonicalb} &=& \left\{\begin{array}{cl}
       1      & \mbox{if $\canonicala=\canonicalb$,} \\
  - (C_\canonicala, C_\canonicalb)
              & \mbox{if $\canonicala<\canonicalb$,} \\
       0      & \mbox{otherwise}
   \end{array}\right. \label{eq:Picard-Lefschetz_Stokes}
\end{eqnarray}
follows.

\section{Intersection numbers of vanishing cycles}
 \label{sc:vanishing_cycles}

\begin{figure}[H]
\centering
\psfrag{v0}{$v_0$}
\psfrag{v1}{$v_1$}
\psfrag{v2}{$v_2$}
\psfrag{v3}{$v_3$}
\psfrag{v4}{$v_4$}
\psfrag{v5}{$v_5$}
\psfrag{v6}{$v_6$}
\psfrag{v7}{$v_7$}
\psfrag{v8}{$v_8$}
\psfrag{v9}{$v_9$}
\psfrag{vequalssomething}{$\begin{array}{rcl}
 v_0 & = & (   0   ,   0   ), \\
 v_1 & = & (   1   ,   0   ), \\
 v_2 & = & (   0   ,   1   ), \\
 v_3 & = & (  -1   ,  -1   ), \\
 v_4 & = & (  -1   ,   0   ), \\
 v_5 & = & (   0   ,  -1   ), \\
 v_6 & = & (  -1   ,   1   ), \\
 v_7 & = & (   1   ,  -1   ), \\
 v_8 & = & (  -1   ,   2   ), \\
 v_9 & = & (   2   ,  -1   ).
	     \end{array}$}
\includegraphics{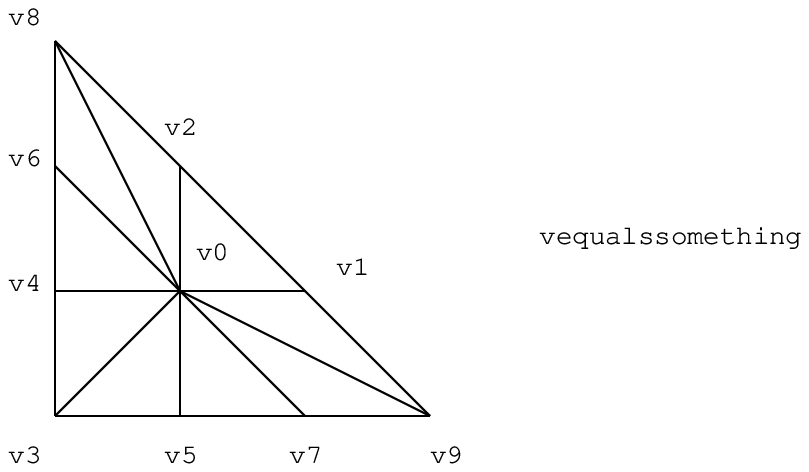}
\caption{the toric data of $\blsix$}
\label{fg:blsix}
\end{figure}

Now we consider the toric variety $\blsix$
which is the projective plane
blown-up at 6 points including infinitely-near points.
Figure \ref{fg:blsix} shows the generators
of the one-dimensional cones
of the fan
defining $\blsix$.

$\blsix$ is deformation-equivalent
to cubic surfaces, and
$-K_\blsix$ is not ample
but nef.
The mirror for $\blsix$ is defined
by (\ref{eq:toric_mirror}),
which we specialize to
$\wa_3 = \wa_8 = \wa_9 =1$ and
$\wa_i = 0$ for $i \neq 3,8,9$,
to obtain
\begin{equation} \label{eq:bl6_mirror}
 W(x,y) = \frac{x^2}{y} + \frac{y^2}{x} + \frac{1}{x y},
\end{equation}
which is generic in the sense
that all the critical points
are non-degenerate.
In this case,
$W^{-1}(0)$ is a Fermat curve in $(\bCx)^2$:
$$
 W^{-1}(0) = \{(x,y) \in (\bCx)^2 \suchthat x^3 + y^3 + 1 = 0 \}.
$$

$W$ has nine critical points
$\{p_{ij}\}_{i,j=0}^2$,
$p_{ij}= (\omega^{i+j-1}, \omega^{2 j})$,
where we have fixed
a primitive cubic root
$\omega$ of unity.
Let $C_{ij}$ be the vanishing cycle
in $W^{-1}(0)$
which vanishes at the critical point $p_{ij}$
along the straight line
from $0$
to the critical value $3 \omega^{1-i}$.
Then
$
(
C_{00},
C_{01},
C_{02},
C_{10},
C_{11},
C_{12},
C_{20},
C_{21},
C_{22}
)
$
is a distinguished basis
of vanishing cycles.

\begin{figure}
 \centering
 \psfrag{C00}{$C_{00}$}
 \psfrag{C01}{$C_{01}$}
 \psfrag{C02}{$C_{02}$}
 \psfrag{C10}{$C_{10}$}
 \psfrag{C11}{$C_{11}$}
 \psfrag{C12}{$C_{12}$}
 \psfrag{C20}{$C_{20}$}
 \psfrag{C21}{$C_{21}$}
 \psfrag{C22}{$C_{22}$}
 \includegraphics{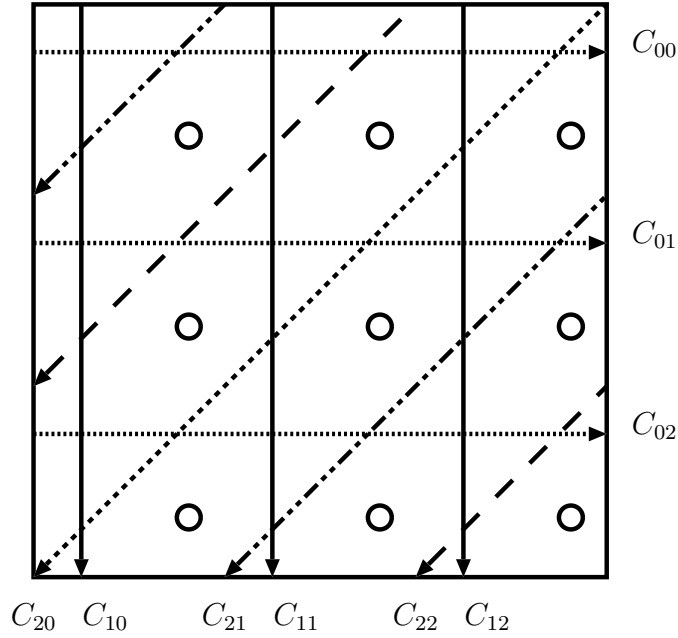}
 \caption{Vanishing cycles in $W^{-1}(0)$}
 \label{fg:vanishing_cycles}
\end{figure}


Figure \ref{fg:vanishing_cycles}
shows these vanishing cycles
with orientations
chosen by hand.
Opposite sides of the square is identified
to form a torus,
and the circles denote nine points
which are missing in $W^{-1}(0)$
since it is not compact.
Note that
the quotient
of $W^{-1}(0)$
by the action of $\bZ / (3 \bZ)$
generated by
$$
 (x, y) \mapsto (\omega x, \omega^2 y).
$$
is
$$
 \left\{
  (X, Y) \in (\bCx)^2 \left| X + Y + \frac{1}{X Y} = 0 \right.
 \right\}
$$
appearing in
Seidel's work
\cite{Seidel_VC2}.
This action is given by the translation
in the direction of the diagonal line
from the lower-left corner
to the upper-right corner
of the square
in Figure \ref{fg:vanishing_cycles}
by one-third of the length
of the diagonal line,
and the quotient
by this action is exactly
Figure 2 in \cite{Seidel_VC2}.

Intersection numbers of $C_{ij}$'s
are
$$
 (C_{ij}, C_{ik}) = 0 \quad \text{if $j \neq k$},
$$
and
$$
 (C_{ij}, C_{kl}) = -1 \quad \text{if $i < k$}.
$$

\section{Derived category of coherent sheaves}
 \label{sc:derived_category}
Let $\bZ / 3 \bZ$ act on a three-dimensional vector space
$V \cong \bC^3$ by
$$
 \bZ / 3 \bZ \ni [1] : V \ni (x:y:z) \mapsto
  (\omega x : \omega^2 y : z) \in V
$$
where $\omega$ is
a primitive cubic root
of unity.
This action defines an action of $\bZ / 3 \bZ$
on $\bP^2 = \bP(V)$,
and the toric variety
$\blsix$ defined in Section \ref{sc:vanishing_cycles}
is the minimal resolution
of the quotient
$\bP(V)/(\bZ / 3 \bZ)$.
By Beilinson \cite{Beilinson},
$(\scE_0, \scE_1, \scE_2)
  = (\scO_{\bP^2}(-1), \Omega_{\bP^2}(1), \scO_{\bP^2})$
is an exceptional collection generating $D^b \coh \bP^2$,
and by Kapranov-Vasserot \cite{Kapranov-Vasserot},
we have
$$
 D^b \coh \blsix \cong D^b \coh^{\bZ/3\bZ}\bP^2.
$$
Here, $\scO_{\bP^2}(-1)$ is the tautological sheaf,
$\Omega_{\bP^2}(1)$ is the cotangent sheaf
tensored with the dual of $\scO_{\bP^2}(-1)$,
and
$\coh^{\bZ/3\bZ}\bP^2$
is the category of $\bZ / 3 \bZ$-equivariant
coherent sheaves on $\bP^2$.
Now, let
$\scE_{ij}$ be the object
of $D^b \coh \blsix$
corresponding to
the $\bZ / 3 \bZ$-equivariant
coherent sheaf
$\scE_i \otimes \rho_j$ on $\bP^2$
by the above equivalence.
Here, $\rho_j : \bZ / 3 \bZ \ni [1] \mapsto \omega^j \in \bCx$
is a one-dimensional representation of $\bZ / 3 \bZ$.
Then
$(
\scE_{00},
\scE_{01},
\scE_{02},
\scE_{10},
\scE_{11},
\scE_{12},
\scE_{20},
\scE_{21},
\scE_{22}
)$
is an exceptional collection
generating $D^b \coh \blsix$.
The $\Ext$-groups between them
can be calculated as follows:
First,
we can use the exact sequence
$$
0 \rightarrow \Omega_{\bP^2}(1) \rightarrow V^\vee \otimes \scO_{\bP^2}
 \rightarrow \scO_{\bP^2}(1) \rightarrow 0
$$
to compute the $\Ext$-groups between $\scE_i$'s:
\begin{eqnarray*}
\RHom(\scE_1, \scE_2)
 & = & \RHom( \scO_{\bP^2}(-1),
         \{ V^\vee \otimes \scO_{\bP^2}
                 \rightarrow \scO_{\bP^2}(1) \} ) \\
 & = & \RGamma( \scO_{\bP^2}(1) \otimes
         \{ V^\vee \otimes \scO_{\bP^2}
                 \rightarrow \scO_{\bP^2}(1) \} ) \\
 & = & \RGamma( \{ V^\vee \otimes \scO_{\bP^2}(1)
                 \rightarrow \scO_{\bP^2}(2) \} ) \\
 & = & \{V^\vee \otimes V^\vee \rightarrow \Sym^2 V^\vee \} \\
 & = & \wedge^2 V^\vee, \\
\RHom(\scE_1, \scE_3)
 &=& \RHom( \scO_{\bP^2}(-1), \scO_{\bP^2} ) \\
 &=& \RGamma( \scO_{\bP^2}(1)) \\
 &=& V^\vee, \\
\RHom(\scE_2, \scE_3)
 & = & \RHom( \{ V^\vee \otimes \scO_{\bP^2}
                 \rightarrow \scO_{\bP^2}(1) \},
               \scO_{\bP^2} ) \\
 & = & \RGamma( \{ \scO_{\bP^2}(-1) \rightarrow
                  V \otimes \scO_{\bP^2} \} ) \\
 & = & V.
\end{eqnarray*}
Here, $\RHom$ and $\RGamma$ denote
the right derived functor of $\Hom$ and $\Gamma$
(taking global sections)
respectively.
Then the $\Ext$-groups between $\scE_{ij}$'s are given by
\begin{eqnarray*}
 \Ext^k (\scE_{ij}, \scE_{lm})
  &=& \Ext^k (\scE_{i} \otimes \rho_j, \scE_{l} \otimes \rho_m) \\
  &=& (\Ext^k (\scE_{i}, \scE_{l})
   \otimes \rho_j^\vee \otimes \rho_m)^{\bZ / 3 \bZ},
\end{eqnarray*}
where $\bullet^{\bZ / 3 \bZ}$ denotes
the $\bZ / 3 \bZ$-invariant part.
Therefore,
the dimension of
$\Ext^k(\scE_{ij}, \scE_{lm})$
is non-zero if and only if
$$
 i=l, j=m \ \text{and} \ k=0,
$$
or
$$
 j<l \ \text{and} \ k=0,
$$
and the dimensions in all these cases
are one.
Therefore,
we have
\begin{eqnarray*}
 -( C_{i j}, C_{l m})
 = \sum_{k=0}^2 (-1)^k \dim \Ext^k (\scE_{i j}, \scE_{l m})
 \quad \text{if $i \leq l$ and $(i,j) \neq (l,m)$},
\end{eqnarray*}
which,
combined with
(\ref{eq:Picard-Lefschetz_Stokes}),
proves Conjecture \ref{conj:stokes}
for $Y$.
\bibliographystyle{plain}
\bibliography{bibs}

Research Institute for Mathematical Sciences,
Kyoto University,
Oiwake-cho,
Kitashirakawa,
Sakyo-ku,
Kyoto,
606-8502,
Japan.

{\em e-mail address}\ : \  kazushi@kurims.kyoto-u.ac.jp

\end{document}